# Discrete random walk with geometric absorption

**Theo van Uem**

School of Technology, Amsterdam University of Applied Sciences,
Weesperzijde 190, 1097 DZ Amsterdam, The Netherlands.

Email: t.j.van.uem@hva.nl

**Abstract**

We consider a discrete random walk (RW) in n dimensions . The RW is adapted with a geometric absorption process: at any discrete time there is a constant probability that absorption occurs in the current state. To model the RW with geometric absorption we use the concept of a multiple function barrier (MFB). In a MFB there is a modification of the original RW: each transition probability in the original RW is multiplied by $\beta$ and there is an additional probability (1-$\beta$) of absorption, where 0<$\beta$<1. We study three cases: one-dimensional simple asymmetric RW , n-dimensional simple symmetric RW (n>1) and a two level RW.

## 1. Introduction

Discrete random walks are studied in a number of standard books, see e.g. Feller [1] and Spitzer [2]. McCrea and Whipple [3] study simple symmetric random paths in two and three dimensions, starting in a rectangular lattice on the integers with absorbing barriers on the boundaries. After taking limits they obtain probabilities of absorption in two and three dimensional lattices. Bachelor and Henry [4,5] use the McCrea-Whipple techniques and find the exact solution for random walks in the triangular lattice with absorbing boundaries and for random walks on finite lattice tubes . We extend the results of Mc Crea and Whipple in two directions: n dimensions and adaption of a geometric absorption process. We consider a discrete random walk (RW) in n dimensions . The RW is adapted with a geometric absorption process: at any discrete time there is a constant probability that absorption occurs in the current state. To model the RW with geometric absorption we use the concept of a multiple function barrier (MFB). In a MFB there is a modification of the original RW: each transition probability in the original RW is multiplied by α and there is an additional probability (1-α) of absorption, where 0<α<1. We study three cases: one-dimensional simple asymmetric RW , n-dimensional simple symmetric RW (n>1) and a two level RW.

## 2. One dimensional simple asymmetric RW with geometric absorption

The original RW starts in 0 and has simple transitions in each state: probability p for one step forward and probability q for one step backward (p+q=1). The adapted RW has probabilities pα (one step forward), qα (one step backward) and (1-α) for absorption in the current state. In this case we have β=α.

The expected number of visits to state j, when starting in 0, is given by: $X_j = \sum_{k=0}^{\infty} p_{0j}^{(k)}$

We have: $\qquad X_n = \delta(n,0) + p\alpha X_{n-1} + q\alpha X_{n+1} \quad (n \in Z)$

**Theorem 1**

The set of difference equations:
$$X_n = \delta(n,0) + p\alpha X_{n-1} + q\alpha X_{n+1} \quad (n \in \mathbb{Z}) \quad (pq > 0,\ p+q = 1,\ 0 < \alpha < 1) \qquad (1)$$

has solutions:
$$X_n = \begin{cases} C_1 \xi_1^n & (n \leq 0) \\ C_2 \xi_2^n & (n \geq 0) \end{cases} \qquad (2)$$

where:
$$\xi_1 = \frac{1 + \sqrt{1 - 4pq\alpha^2}}{2q\alpha} > 1$$

$$0 < \xi_2 = \frac{1 - \sqrt{1 - 4pq\alpha^2}}{2q\alpha} < 1$$

$$C_1 = C_2 = [1 - 4pq\alpha^2]^{-\frac{1}{2}}$$

**Proof**

General solution of homogeneous part of (1) is:
$$X_n = C_1 \xi_1^n + C_2 \xi_2^n \quad (n \in \mathbb{Z}),\ \text{where}\ q\alpha \xi^2 - \xi + p\alpha = 0$$

$X_n$ is finite ($\alpha > 0$), so we have (2).

$C_1$ and $C_2$ are determined by: using (2), $X_0 = X_0$ and $X_0 = 1 + p\alpha X_{-1} + q\alpha X_1$ □

Absorption probabilities can be obtained by $P(absorption\ in\ state\ n) = (1 - \alpha)X_n$

## 3. n-dimensional simple symmetric RW with geometric absorption (n>1)

We again use $X_j = \sum_{k=0}^{\infty} p_{0j}^{(k)}$, where 0 and j are now vectors in a n-dimensional lattice on the integers.

In each point there is probability $\alpha$ to move one step in any of the 2n directions and probability (1-2n$\alpha$) of absorption ($0 < \alpha < \frac{1}{2n}$). In this case we have $\beta$=2n$\alpha$.

We have:
$$X_u = \delta(u,0) + \sum_{k=1}^{n} \alpha \{X_{u-[k]} + X_{u+[k]}\}$$

where [k] is the n-dimensional vector of length one in positive direction k : $[k] = (u_1, u_2, \ldots\ldots u_n)$, $u_i = \delta(k, 0)$, $i = 1, 2, \ldots, n$.

**Theorem 2**

The set of difference equations:
$$X_u = \delta(u,0) + \sum_{k=1}^{n} \alpha \{X_{u-[k]} + X_{u+[k]}\} \qquad (3)$$

where 0 and u are n-dimensional vectors in an integer valued lattice and [k] is the n-dimensional vector of length one in the positive direction k, has solution:

$$X_u = \frac{1}{2\alpha\pi^{n-1}} \int_0^\pi \cdots \int_0^\pi \prod_{i=1}^{n-1} cos u_i \omega_i \frac{e^{-|u_n|}}{sinh\omega_n} d\omega_1 \ldots\ldots\ldots d\omega_{n-1} \quad (4)$$

where
$$\sum_{k=1}^{n-1} \cos\omega_i + \cosh\omega_n = \frac{1}{2\alpha} \quad (5)$$

Proof
That (4) under condition (5) is a solution of (3) can easily be verified by direct substitution and using elementary goniometrical identities. □

Absorption probabilities can now be obtained by $P(absorption\ in\ state\ u) = (1 - 2n\alpha)X_u$

## 4. A two level RW with geometric absorption

Both levels consist of the set of integers. We start at level zero in 0. At any level we can jump with probability $\alpha$ to the direct neighbors at the same level or to the direct neighbor on the other level:

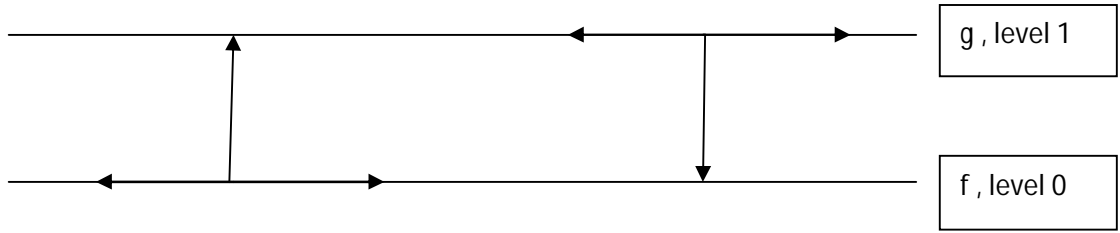

In each point of both levels there is a probability $(1 - 3\alpha)$ of absorption. We now have β=3α.
We describe the behavior at level 0 with function f and at level 1 with function g:

$$f_j = \sum_{k=0}^\infty p_{0j}^{(k)} \text{ (0 and j both on level 0)};\quad g_j = \sum_{k=0}^\infty p_{0j}^{(k)} \text{ (j on level 1, 0 on level 0)}.$$

We readily get:
$$f_n = \delta(n,0) + \alpha(f_{n-1} + f_{n+1} + g_n) \quad (n \in Z) \quad (6)$$
$$g_n = \alpha(g_{n-1} + g_{n+1} + f_n) \quad (n \in Z) \quad (7)$$

Combining this equations yields:
$$f_{n+2} - \frac{2}{\alpha}f_{n+1} + \left(1 + \frac{1}{\alpha^2}\right)f_n - \frac{2}{\alpha}f_{n-1} + f_{n-2} = 0$$

The characteristic equation can be written as:
$$\left(\mu^2 + \left(1 - \frac{1}{\alpha}\right)\mu + 1\right)\left(\mu^2 + \left(-1 - \frac{1}{\alpha}\right)\mu + 1\right) = 0$$

with solutions:
$$\mu_{1,2} = \frac{1 - \alpha \pm \sqrt{(1+\alpha)(1-3\alpha)}}{2\alpha}$$

$$\mu_{3,4} = \frac{1 + \alpha \pm \sqrt{(1-\alpha)(1+3\alpha)}}{2\alpha}$$

Because of $\mu_1$ and $\mu_3$ both greater than 1 and using the symmetry, we have:
$$f_n = \begin{cases} b\mu_2^n + d\mu_4^n & (n \geq 0) \\ b\mu_2^{-n} + d\mu_4^{-n} & (n \leq 0) \end{cases}$$

Using (6) we get:

$$g_n = \begin{cases} b\mu_2^n - d\mu_4^n & (n \geq 0) \\ b\mu_2^{-n} - d\mu_4^{-n} & (n \leq 0) \end{cases}$$

At last we use (6) and (7) both with n=0, which gives us:

$$b = \frac{1}{2\sqrt{(1+\alpha)(1-3\alpha)}}$$

$$d = \frac{1}{2\sqrt{(1-\alpha)(1+3\alpha)}}$$

Absorption probabilities can now be obtained by $P(absorption\ in\ state\ n\ on\ level\ 0) = (1 - 3\alpha)f_n$ and $P(absorption\ in\ state\ n\ on\ level\ 1) = (1 - 3\alpha)g_n$

### 5. Uniqueness of solutions

We have found solutions of different problems, but are they unique?

In all cases we have transient random walk (positive constant probability of absorption at any time).

We prove the uniqueness in the n-dimensional case. The two level case can be handled in a similar way.

Let $u = (u_1, u_2, \ldots, u_n)$. We have: $X_u \to 0$ if $\sum_{i=1}^{n} u_i^2 \to \infty$, so we get:

$\forall \epsilon > 0 \ \exists R_\epsilon : 0 \leq X_u < \epsilon$ if $\sum_{i=1}^{n} u_i^2 > R_\epsilon^2$

Let M be the supremum of X on the n-dimensional lattice. Suppose M>0.

$\forall \epsilon > 0, M > 2\epsilon \quad \exists u^M = (u_1^M, u_2^M, \ldots, u_n^M) : X_{u^M} > M - \epsilon > \epsilon > 0$

Combining the two last results we find that M can only be found within the n-sphere on the integers with radius $R_\epsilon$: M is the maximum on that set.

Because of $\alpha > 0$ it follows now from (3) that $X_{[n]+u^M} = X_{u^M}$

Continuing this way gives:

$\forall k \epsilon N : X_{k[n]+u^M} > \epsilon$ and $\sum_{i=1}^{n-1}(u_i^M)^2 + (u_n^M + k)^2 < R_\epsilon^2$ where k[n] is k times the unit vector in positive direction n. Contradiction.

Conclusion: The homogeneous part of (3) has only the solution 0.